# $L^p$-differentiability: A unifying theme in characterizing spaces of weakly differentiable functions

Daniel Spector*

Department of Mathematics, Technion - Israel Institute of Technology

June 3, 2014

**Abstract**

In this note, we introduce a variant of the notion of $L^p$-differentiability introduced by Calderón and Zygmund and prove that functions in a Sobolev space possess this type of $L^p$-derivative. In contrast to the classical notion, we show that our formulation of the condition is not only a property of the Sobolev functions, but in fact characterizes them. Our approach unifies some characterizations of Sobolev spaces due to Swanson using Calderón-Zygmund classes with others due to Bourgain, Brezis, and Mironescu using nonlocal functionals and the author and Mengesha using nonlocal gradients. That any two characterizations of Sobolev spaces are related is not surprising, however, one consequence of our analysis is a simple condition for determining whether a function of bounded variation is in a Sobolev space.

## 1  Introduction and Main Results

$L^p$-differentiability was introduced by Calderón and Zygmund in their study of the local properties of solutions of elliptic differential equations [4, 5]. It is a natural extension of classical differentiability in that it relaxes the requirement of the existence of a locally linear map in a uniform sense to its existence in an averaged sense. As the Sobolev spaces arise readily in the study of partial differential equations, it is not surprising that Sobolev functions possess an $L^p$-derivative. The following theorem asserting this fact was proven[1] by Calderón and Zygmund in [5].

**Theorem 1.1 (Calderón and Zygmund)** *Suppose $1 \leq p < \infty$ and $f \in W^{1,p}(\mathbb{R}^N)$. Then*

$$\left(\frac{1}{\epsilon^{pq}} \fint_{B(0,\epsilon)} |f(x+h) - f(x) - \nabla f(x) \cdot h|^{pq} \, dh\right)^{\frac{1}{pq}} \to 0 \qquad (1.1)$$

---

*dspector@tx.technion.ac.il

[1] For a modern reference, we refer to the monograph of Evans and Gariepy [7].

for $\mathcal{L}^N$ almost every $x \in \mathbb{R}^N$, where $1 \leq q \leq \frac{N}{N-p}$ if $1 \leq p < N$, $1 \leq q < \infty$ if $p = N$, and $1 \leq q \leq \infty$ if $p > N$ (Here, $q = +\infty$ is understood to be the $L_h^\infty(B(0,\epsilon))$ norm applied to the integrand).

In the language of Calderón and Zymund, Theorem 1.1 states that $f \in W^{1,p}(\mathbb{R}^N)$ implies $f \in t^{1,pq}(x)$ (which is essentially defined by the condition (1.1)) for almost every $x \in \mathbb{R}^N$. The converse is false, readily seen through the "improved" $L^{pq}$-differentiability of functions in the Sobolev space $W^{1,p}(\mathbb{R}^N)$. A natural question is then whether one can characterize the Sobolev spaces in the spirit of the condition (1.1). Several results have been given in this direction using the class[2] $T^{1,p}(x)$, also introduced by Calderón and Zygmund. We mention, for instance, a sort of converse to Theorem 1.1 due to Bagby and Ziemer [2], as well as work by Swanson [10, 11], which additionally assumes either integrability or an infinitesimal property (and is actually related to our viewpoint). Here we will begin by examining the condition (1.1) from a different perspective.

Our approach begins with the observation that the key ingredients to Theorem 1.1 are the Sobolev embedding theorem and the Lebesgue differentiation theorem for $L^p$ functions. As in the statement of Theorem 1.1, the latter is typically stated as a pointwise almost everywhere convergence result. However, a variant of the Lebesgue differentiation theorem, which has been pointed out to us by Michael Cwickel, is the following $L^p(\mathbb{R}^N)$ convergence result. If $u \in L^p(\mathbb{R}^N)$ and we define

$$h_\epsilon(x) := \left( \fint_{B(x,\epsilon)} |u(x) - u(y)|^p \, dy \right)^{\frac{1}{p}},$$

then one has the convergence $h_\epsilon \to 0$ in $L^p(\mathbb{R}^N)$ as $\epsilon \to 0$.

From this perspective, it would be natural to expect an analogous $L^p$-type convergence to be true for functions in the Sobolev space $W^{1,p}(\mathbb{R}^N)$. The first result we announce in this note is the following theorem to this effect.

**Theorem 1.2** Let $1 \leq p < \infty$ and $f \in W^{1,p}(\mathbb{R}^N)$. Then

$$\lim_{\epsilon \to 0} \int_{\mathbb{R}^N} \left( \fint_{B(0,\epsilon)} \frac{|f(x+h) - f(x) - \nabla f(x) \cdot h|^{pq}}{|h|^{pq}} \, dh \right)^{\frac{1}{q}} dx = 0, \qquad (1.2)$$

where $1 \leq q \leq \frac{N}{N-p}$ if $1 \leq p < N$, $1 \leq q < \infty$ if $p = N$, and $1 \leq q \leq \infty$ if $p > N$.

Here we have utilized a variant of the integrand in (1.1) in the statement of our theorem, whose pointwise convergence to zero is equivalent to the convergence (1.1), a fact which could be deduced from the subsequent analysis we present in this paper (see also [1][Chapter , Exercise ]).

What is quite surprising, and one of the main results we announce here, is that this property in fact characterizes Sobolev functions. Precisely, we have the following theorem characterizing the Sobolev space $W^{1,p}(\mathbb{R}^N)$ in terms of the $L^p(\mathbb{R}^N)$ convergence (1.2).

---
[2] For precise definitions of $t^{1,p}(x)$ and $T^{1,p}(x)$, we refer to Ziemer [12][Chapter 3, p.132].

**Theorem 1.3** *Let $1 \leq p < +\infty$ and suppose $f \in L^p(\mathbb{R}^N)$. Then $f \in W^{1,p}(\mathbb{R}^N)$ if and only if there exists a $v \in L^p(\mathbb{R}^N; \mathbb{R}^N)$ such that*

$$\lim_{\epsilon \to 0} \int_{\mathbb{R}^N} \left( \fint_{B(0,\epsilon)} \frac{|f(x+h) - f(x) - v(x) \cdot h|^{pq}}{|h|^{pq}} \, dh \right)^{\frac{1}{q}} dx = 0 \qquad (1.3)$$

*for any $1 \leq q \leq \frac{N}{N-p}$ if $1 \leq p < N$, $1 \leq q < \infty$ if $p = N$, and $1 \leq q \leq \infty$ if $p > N$.*

Actually, for $1 < p < +\infty$, we will see in the proof that the finiteness of the limit (1.3) is sufficient to deduce that $f \in W^{1,p}(\mathbb{R}^N)$. Such an assumption is related to the characterization of Swanson in [10], while the convergence to zero is analogous to an assumption in the paper [11] which allows the author to characterize $W^{1,1}(\mathbb{R}^N)$. When $p = 1$, if one only assumes that the limit (1.3) is finite, one instead obtains a characterization of the space of functions of bounded variation in the spirit of the characterization of Bourgain, Brezis, and Mironescu [3], which is the content of the following theorem.

**Theorem 1.4** *Suppose $f \in L^1(\mathbb{R}^N)$. Then $f \in BV(\mathbb{R}^N)$ if and only if there exists a $v \in L^1(\mathbb{R}^N; \mathbb{R}^N)$ such that*

$$\lim_{\epsilon \to 0} \int_{\mathbb{R}^N} \left( \fint_{B(0,\epsilon)} \frac{|f(x+h) - f(x) - v(x) \cdot h|^q}{|h|^q} \, dh \right)^{\frac{1}{q}} dx < +\infty \qquad (1.4)$$

*for any $1 \leq q \leq \frac{N}{N-1}$.*

As a consequence of Theorems 1.4 and 1.3, we have that a function $f \in BV(\mathbb{R}^N)$ is included in $W^{1,1}(\mathbb{R}^N)$ if and only if (1.3) holds for $p = 1$. For a general $f \in BV(\mathbb{R})$, this implies the limit in equation (1.4) for $q = 1$ gives an upper bound for the total mass of the singular part of the measure $Df$.

We will shortly give a proof of Theorems 1.2 and 1.3 in the case $1 \leq p < N$, while the proofs of the regime $p \geq N$ and Theorem 1.4 will be deferred to a later work. First, let us state without proof the following Lemma, which is a variation of a calculation implicit in Evans and Gariepy [7][Chapter 6, p. 231].

**Lemma 1.5** *Suppose $f \in W^{1,p}_{loc}(\mathbb{R}^N)$ for some $1 \leq p < \infty$, and that $1 \leq q \leq \frac{N}{N-p}$ if $1 \leq p < N$. Then there exists a $C = C(p, q, N) > 0$ such that for all $0 < t < 1$*

$$\frac{1}{t^{N+pq}} \int_{B(0,t)} |f(x+h) - f(x) - \nabla f(x)h|^{pq} \, dh$$

$$\leq C \left( \fint_{B(0,t)} |\nabla f(x+h) - \nabla f(x)|^p \, dh \right)^q$$

$$+ C \left( \int_0^1 \fint_{B(0,st)} |\nabla f(x+sz) - \nabla f(x)|^p \, dz ds \right)^q.$$

We now give a proof of Theorem 1.2.

**Proof.** For any $0 < \epsilon < 1$, we expand the integrand on concentric rings

$$\fint_{B(0,\epsilon)} \frac{|f(x+h) - f(x) - \nabla f(x) \cdot h|^{pq}}{|h|^{pq}} \, dh$$

$$= \sum_{i=0}^{\infty} \frac{1}{\epsilon^N |B(0,1)|} \int_{B(0,\frac{\epsilon}{2^i}) \setminus B(0,\frac{\epsilon}{2^{i+1}})} \frac{|f(x+h) - f(x) - \nabla f(x) \cdot h|^{pq}}{|h|^{pq}} \, dh.$$

We now make estimates for $i \in \mathbb{N}$ fixed. We have

$$\frac{1}{\epsilon^N |B(0,1)|} \int_{B(0,\frac{\epsilon}{2^i}) \setminus B(0,\frac{\epsilon}{2^{i+1}})} \frac{|f(x+h) - f(x) - \nabla f(x) \cdot h|^{pq}}{|h|^{pq}} \, dh$$

$$\leq \frac{1}{\epsilon^N |B(0,1)|} \left(\frac{\epsilon}{2^{i+1}}\right)^{-pq} \int_{B(0,\frac{\epsilon}{2^i}) \setminus B(0,\frac{\epsilon}{2^{i+1}})} |f(x+h) - f(x) - \nabla f(x) \cdot h|^{pq} \, dh$$

$$= \frac{2^{pq}}{|B(0,1)|} \frac{1}{2^{iN}} \left(\frac{\epsilon}{2^i}\right)^{-N-pq} \int_{B(0,\frac{\epsilon}{2^i})} |f(x+h) - f(x) - \nabla f(x) \cdot h|^{pq} \, dh.$$

Lemma 1.5 further implies that

$$\left(\frac{\epsilon}{2^i}\right)^{-N-pq} \int_{B(0,\frac{\epsilon}{2^i})} |f(x+h) - f(x) - \nabla f(x) \cdot h|^{pq} \, dh$$

$$\leq C \left(\fint_{B(0,\frac{\epsilon}{2^i})} |\nabla f(x+h) - \nabla f(x)|^p \, dh\right)^q$$

$$+ C \left(\int_0^1 \fint_{B(0,s\frac{\epsilon}{2^i})} |\nabla f(x+sz) - \nabla f(x)|^p \, dzds\right)^q.$$

Therefore, summing in $i$ and applying the basic inequality $(\sum_i |a_i|)^{\frac{1}{q}} \leq \sum_i |a_i|^{\frac{1}{q}}$, we have

$$\left(\fint_{B(0,\epsilon)} \frac{|f(x+h) - f(x) - \nabla f(x)h|^{pq}}{|h|^{pq}} \, dh\right)^{\frac{1}{q}}$$

$$\leq C \sum_{i=0}^{\infty} \left(\frac{1}{2^i}\right)^{N/q} \fint_{B(0,\frac{\epsilon}{2^i})} |\nabla f(x+h) - \nabla f(x)|^p \, dh$$

$$+ C \sum_{i=0}^{\infty} \left(\frac{1}{2^i}\right)^{N/q} \int_0^1 \fint_{B(0,s\frac{\epsilon}{2^i})} |\nabla f(x+sz) - \nabla f(x)|^p \, dzds.$$

Integrating the preceding inequality over $x \in \mathbb{R}^N$ and making use of Tonelli's theorem we obtain

$$\int_{\mathbb{R}^N} \left(\fint_{B(0,\epsilon)} \frac{|f(x+h) - f(x) - \nabla f(x)h|^{pq}}{|h|^{pq}} \, dh\right)^{\frac{1}{q}} dx$$

$$\leq C \sum_{i=0}^{\infty} \left(\frac{1}{2^i}\right)^{N/q} \fint_{B(0,\frac{\epsilon}{2^i})} \int_{\mathbb{R}^N} |\nabla f(x+h) - \nabla f(x)|^p \, dxdh$$

$$+ C \sum_{i=0}^{\infty} \left(\frac{1}{2^i}\right)^{N/q} \int_0^1 \fint_{B(0,s\frac{\epsilon}{2^i})} \int_{\mathbb{R}^N} |\nabla f(x+sz) - \nabla f(x)|^p \, dxdzds.$$

However, if $h, z \in B(0, \frac{\epsilon}{2^i})$ and $s \in (0,1)$ we have

$$\max\left\{\int_{\mathbb{R}^N} |\nabla f(x+h) - \nabla f(x)|^p\, dx, \int_{\mathbb{R}^N} |\nabla f(x+sz) - \nabla f(x)|^p\, dx\right\}$$
$$\leq \sup_{w \in B(x,\epsilon)} \int_{\mathbb{R}^N} |\nabla f(x+w) - \nabla f(x)|^p\, dx,$$

and observe that this bound is independent of $i \in \mathbb{N}$. Thus,

$$\int_{\mathbb{R}^N} \left(\fint_{B(0,\epsilon)} \frac{|f(x+h) - f(x) - \nabla f(x)h|^{pq}}{|h|^{pq}}\, dh\right)^{\frac{1}{q}} dx$$
$$\leq \sup_{w \in B(x,\epsilon)} \int_{\mathbb{R}^N} |\nabla f(x+w) - \nabla f(x)|^p\, dx \left(C \sum_{i=0}^{\infty} \left(\frac{1}{2^i}\right)^{N/q}\right).$$

As the infinite series is summable, the result follows from sending $\epsilon \to 0$ and using continuity of translation in $L^p(\mathbb{R}^N)$. ∎

We conclude with a proof of Theorem 1.3.

**Proof.** As we have shown that $f \in W^{1,p}(\mathbb{R}^N)$ implies the $L^p$-convergence (1.3), it remains to show the converse. We first treat the case $1 < p < +\infty$. Let us therefore suppose that there exists $v \in L^p(\mathbb{R}^N; \mathbb{R}^N)$ such that (1.3) holds. We then estimate

$$\int_{\mathbb{R}^N} \fint_{B(x,\epsilon)} \frac{|f(x+h) - f(x)|^p}{|h|^p}\, dhdx$$
$$\leq 2^{p-1} \int_{\mathbb{R}^N} \fint_{B(x,\epsilon)} \frac{|f(x+h) - f(x) - v(x) \cdot h|^p}{|h|^p}\, dhdx$$
$$+ 2^{p-1} \int_{\mathbb{R}^N} \fint_{B(x,\epsilon)} \left|v(x) \cdot \frac{h}{|h|}\right|^p dhdx.$$

Now our assumption (and we can take $q=1$) is that the first term on the right hand side tends to zero as $\epsilon \to 0$, while the second is bounded by a constant times the $L^p$ norm of $v$. We then have that for a sequence $\epsilon_n \to 0$

$$\limsup_{n \to \infty} \int_{\mathbb{R}^N} \fint_{B(x,\epsilon_n)} \frac{|f(x+h) - f(x)|^p}{|h|^p}\, dhdx < +\infty,$$

and so by the result of Bourgain, Brezis, and Mironescu [3] we conclude that $f \in W^{1,p}(\mathbb{R}^N)$.

For the case $p=1$, a similar argument allows us to deduce that $f \in BV(\mathbb{R}^N)$ (and can be used to demonstrate one direction of Theorem 1.4). It therefore remains to show that $v = Df$ in the sense of distributions and we can conclude $f \in W^{1,1}(\mathbb{R}^N)$. However, we observe that $f \in BV(\mathbb{R}^N)$ implies that the nonlocal gradient

$$\mathcal{G}_\epsilon f(x) := N \fint_{B(x,\epsilon)} \frac{f(x+h) - f(x)}{|h|} \frac{h}{|h|}\, dhdx$$

is well-defined as a Lebesgue integral and that $\mathcal{G}_\epsilon f \in L^1(\mathbb{R}^N)$ (see the paper of the author and Mengesha, [9]). If we could show that $\mathcal{G}_\epsilon f \to v$ in $L^1(\mathbb{R}^N; \mathbb{R}^N)$, we would be finished, since the convergence $\mathcal{G}_\epsilon f \stackrel{*}{\rightharpoonup} Df$ in $\bigl(C_0(\mathbb{R}^N;\mathbb{R}^N)\bigr)'$ would then imply that $v = Df$ in the sense of distributions, which is the desired result. However, we observe that

$$v_i(x) N \fint_{B(0,\epsilon)} \frac{h_i h_j}{|h|^2} = v_i(x) \delta_{ij}$$

and therefore we have

$$(\mathcal{G}_\epsilon f)_i(x) - v_i(x) = N \fint_{B(0,\epsilon)} \frac{f(x+h) - f(x) - v(x) \cdot h}{|h|} \frac{h_i}{|h|}\, dh,$$

and thus we can estimate

$$\int_{\mathbb{R}^N} |(\mathcal{G}_\epsilon f)_i(x) - v_i(x)| \leq N \int_{\mathbb{R}^N} \fint_{B(0,\epsilon)} \frac{|f(x+h) - f(x) - v(x) \cdot h|}{|h|}\, dh,$$

which tends to zero by our assumption and the result is demonstrated. ∎

In forthcoming work we will complete the proof of Theorem 1.4, as well as to discuss several applications and variations of the result in the context of the assumptions of Bourgain, Brezis, and Mironescu [3]. In particular, we will address the use of different approximations of the identity, local convergence results and the case of characterizations for domains $\Omega \subset \mathbb{R}^N$ open, as well as to give a proof of a result related to a claim in the paper [8].

## Acknowledgements

The author would like to thank Michael Cwikel, Rahul Garg, and Itai Shafrir for their insightful discussions in the process of preparing this result for publication. The author is supported in part by a Technion Fellowship, and is deeply appreciative of Liang-Kuang King and Leh-Ping Yang for their financial support during the preparation of this paper.